\input ./style/arxiv-general.cfg
\documentclass[aop,MSNbibl,dvips]{arximspdf}
\makeatletter
   \@ifpackageloaded{graphicx}{}{\usepackage{graphicx}}
\makeatother
\usepackage{mathbh}

\doi{10.1214/14-AOP990}% Updated by VTEXPTS2LaTeX.exe, 23.12.2014 11:01
\volume{44}
\issue{2}
\pubyear{2016}
\firstpage{867}
\lastpage{882}
\docsubty{FLA}

\makeatletter
\newcommand{\rrvert}{\vert}
\newcommand{\rrVert}{\Vert}
\newcommand{\llvert}{\vert}
\newcommand{\llVert}{\Vert}
\renewcommand{\mid}{|}
\newtheorem{teo}{Theorem}%[section]
\newtheorem{prop}[teo]{Proposition}
\newtheorem{cor}[teo]{Corollary}
\newtheorem{lem}[teo]{Lemma}
\newproclaim{remark}[teo]{Remark}

\newcommand{\real}{{\mathbb R}}
\newcommand{\nat}{{\mathbb N}}
\newcommand{\com}{{\mathbb C}}
\newcommand{\un}{{\mathbh{1}}}
\newcommand{\FF}{{\mathbb F}}

\newcommand{\B}{{\mathcal B}}
\newcommand{\E}{{\mathcal E}}
\newcommand{\M}{{\mathcal M}}
\newcommand{\N}{{\mathcal N}}

\newcommand{\g}{\gamma}
\renewcommand{\d}{\delta}
\renewcommand{\t}{\tau}
\newcommand{\f}{\varphi}
\renewcommand{\l}{\lambda}
\renewcommand{\O}{\Omega}
\newcommand{\s}{\sigma}
\newcommand{\el}{\ell}
\newcommand{\ma}{{\mathbb M}}
\makeatother

\begin{document}
\begin{frontmatter}

\title{A noncommutative martingale convexity inequality\thanksref{T1}}
\runtitle{A noncommutative martingale convexity inequality}

\begin{aug}
% Corresponding author: Quanhua Xu - qxu@univ-fcomte.fr% Updated by
%VTEXPTS2LaTeX.exe, 23.12.2014 11:01
\author[A]{\fnms{\'Eric}~\snm{Ricard}\thanksref{M1}\ead[label=e1]{eric.ricard@unicaen.fr}}
\and
\author[B]{\fnms{Quanhua}~\snm{Xu}\thanksref{M2,M3}\corref{}\ead[label=e2]{qxu@univ-fcomte.fr}}
\runauthor{E. Ricard and Q. Xu}
\affiliation{Universit{\'e} de Caen Basse-Normandie\thanksmark{M1},\\
Wuhan University\thanksmark{M2} and
Universit{\'e} de Franche-Comt{\'e}\thanksmark{M3}}
%\dedicated{}
\address[A]{Laboratoire de Math{\'e}matiques Nicolas Oresme\\
Universit{\'e} de Caen Basse-Normandie\\
14032 Caen Cedex\\
France\\
\printead{e1}}
\address[B]{School of Mathematics and Statistics\\
Wuhan University\\
Wuhan 430072\\
China\\
and\\
Laboratoire de Math{\'e}matiques\\
Universit{\'e} de Franche-Comt{\'e}\\
25030 Besan\c{c}on Cedex\\
France\\
\printead{e2}}
\end{aug}
\thankstext{T1}{Supported by ANR-2011-BS01-008-01 and NSFC Grant No. 11271292.}

% HISTORY:
%
\received{\smonth{5} \syear{2014}}% Updated by VTEXPTS2LaTeX.exe,
%23.12.2014 11:01
%
\revised{\smonth{11} \syear{2014}}% Updated by VTEXPTS2LaTeX.exe,
%23.12.2014 11:01

% ABSTRACT
%
\begin{abstract}
Let $\mathcal M$ be a von Neumann algebra equipped with a faithful
semifinite normal weight
$\phi$ and $\mathcal N$ be a von Neumann subalgebra
of $\mathcal M$ such that the restriction of $\phi$ to $\mathcal N$ is
semifinite and such that $\mathcal N$ is invariant
by the modular group of $\phi$. Let $\mathcal E$ be the weight
preserving conditional expectation
from $\mathcal M$ onto $\mathcal N$. We prove the following inequality:
\[
\| x\|_{p}^2\ge\bigl\|\mathcal E(x)\bigr\|_p^2+(p-1)\bigl\| x-\mathcal E(x)\bigr\|
_p^2,\qquad x\in L_p(\mathcal M), 1<p\le2,
\]
which extends the celebrated Ball--Carlen--Lieb convexity inequality.
As an application we show that there exists $\varepsilon_0>0$
such that for any free group $\mathbb F_n$ and any $q\ge4-\varepsilon_0$,
\[
\| P_t\|_{2\to q}\le1\quad\Leftrightarrow\quad t\ge\log{
\sqrt{q-1}},
\]
where $(P_t)$ is the Poisson semigroup defined by the natural length
function of $\mathbb F_n$.
\end{abstract}

% KEYWORDS
\end{frontmatter}

%s1 #&#
\section{Introduction}

%%%%%%%%%%%%%%%%%%%%%%%%%%%%%%%%%%%%%%%%%%%%%%%%%%%%%%
%%%%%%%%%%%%%%%%%%%%%%%%%%%%%%%%%%%%%%%%%%%%%%%%%%%%%%

Let $\M$ be a von Neumann algebra equipped with a faithful semifinite
normal weight
$\phi$. The associated noncommutative $L_p$-spaces will be simply
denoted by $L_p(\M)$. We refer to \cite{px} for information on
noncommutative integration.
Recall that if $\N$ is a von Neumann subalgebra
of $\M$ such that the restriction of $\phi$ to $\N$ is semifinite
and such that $\N$ is $\s^\phi$-invariant [i.e.,
$\sigma_t^\phi(\N)=\N$ for all $t\in\real$], then there
exists a unique $\phi$-preserving conditional expectation $\E$ from
$\M$ onto $\N$ such that
\[
\E(axb)=a\E(x)b,\qquad a, b\in\N, x\in\M.
\]
Here $\sigma^\phi$ denotes the modular group of $\phi$. Moreover,
$\E$ extends
to a contractive projection from $L_p(\M)$ onto $L_p(\N)$ for any
$1\le p<\infty$. Below is our main result.

%th1 #&#
\begin{teo}\label{conv-th}
Let $\M$, $\N$ and $\E$ be as above. If $1<p\le2$, then
%
%e1 #&#
\begin{equation}
\label{conv} \llVert x\rrVert_{p}^2\ge\bigl\llVert\E(x)
\bigr\rrVert_p^2+(p-1)\bigl\llVert x-\E(x)\bigr\rrVert
_p^2,\qquad x\in L_p(\M).
\end{equation}
If $2<p<\infty$, the inequality is reversed.
\end{teo}

Inequality (\ref{conv}) is a martingale convexity
inequality. It is closely related to the celebrated convexity
inequality of Ball, Carlen and Lieb \cite{bcl} for the Schatten
classes~$S_p$. Namely, for $1<p\le2$, we have
%
%e2 #&#
\begin{equation}
\label{bcl} \llVert x+y\rrVert^2_p+ \llVert x-y\rrVert
_p^2\geq2 \llVert x\rrVert_p^2
+2(p-1) \llVert y\rrVert_p^2,\qquad x, y\in
S_p.
\end{equation}
In fact, it is easy to see that (\ref{bcl}) is a
special case of (\ref{conv}) by considering $\M=B(\el_2)\oplus
B(\el_2)$. Conversely, the validity of~(\ref{bcl}) for any
noncommutative $L_p$-spaces implies~(\ref{conv}). Indeed, we will
deduce (\ref{conv}) from the following:

%th2 #&#
\begin{teo}\label{BCL-th}
Let $\M$ be any von Neumann algebra. If $1<p\le2$, then
%
%e3 #&#
\begin{equation}
\label{BCL} \llVert x+y\rrVert^2_p+ \llVert x-y\rrVert
_p^2\geq2 \llVert x\rrVert_p^2
+2(p-1) \llVert y\rrVert_p^2,\qquad x, y\in
L_p(\M).
\end{equation}
If $2<p<\infty$, the inequality is reversed.
\end{teo}

What is new and remarkable in (\ref{bcl}) or (\ref{BCL}) is the fact
that $(p-1)$ is the best constant. In fact, if one allows a constant
depending on $p$ in place of $(p-1)$, then~(\ref{BCL}) is equivalent
to the well-known results on the 2-uniform convexity of $L_p(\M)$. We
refer to \cite{bcl} for more discussion on this point. The optimality
of the constant $(p-1)$ has important applications to
hypercontractivity in the noncommutative case. It is the key to the
solution of Gross's longstanding open problem about the optimal
hypercontractivity for Fermi fields by Carlen and Lieb \cite{cl}. It
plays the same role in \cite{b} and \cite{jpppr}. Note that the
optimality of $(p-1)$ in (\ref{BCL}) implies $(p-1)$ is also the best
constant in (\ref{conv}). It seems that (\ref{conv}) with this best
constant is new even in the commutative case.

Clearly, (\ref{bcl}) implies (\ref{BCL}) for injective $\M$ (or more
generally, QWEP $\M$) since then $L_p(\M)$ is finitely representable
in $S_p$. The proof of (\ref{bcl}) in \cite{bcl} goes through a
differentiation argument for the function $t\mapsto
\llVert x+ty\rrVert _p^p$ with self-adjoint $x$ and~$y$. It seems
difficult to directly
extend their argument to finite von
Neumann algebras. The subtle point is the fact that to be
able to differentiate the above function, one needs the invertibility
of $x+ty$ for all $t\in[0, 1]$ except possibly countably many of
them. This
invertibility is easily achieved in the matrix algebra case, that is,
for $\M=\ma_n$, the algebra of $n\times n$-matrices. Instead, we will
use a pseudo-differentiation argument which is much less rigid than
that of \cite{bcl}. The main novelty in our argument can be simply
explained as follows. We first cut the operator $x+ty$ by its spectral
projections in order to reduce the general case to the invertible
one; to do so we need $x+ty$ to be of full support for all
$t\in[0, 1]$. We then get this full support property for all $t$ by
adding to $x+ty$ an independent operator with diffuse spectral
measure. Note that by standard perturbation argument, it is easy to
insure the full support (or even the invertibility) of $x+ty$ for one
$t$.

An iteration of Theorem~\ref{conv-th} immediately implies the
following inequality
on noncommutative martingales.

%co3 #&#
\begin{cor}
Let $(\M_n)_{n\ge0}$ be an increasing sequence of von Neumann
subalgebras of $\M$
with w*-dense union in $\M$. Assume that each $\M_n$ is \mbox{$\s^\phi$-}invariant,
and $\phi\mid_{\M_n}$ is semifinite. Let $\E_n$ be the conditional
expectation
with respect to $\M_n$. Then for $1<p\le2$,
\[
\llVert x\rrVert_p^2\ge\bigl\llVert
\E_0(x)\bigr\rrVert_p^2+(p-1)\sum
_{n\ge1}\bigl\llVert\E_n(x)-\E_{n-1}(x)
\bigr\rrVert_p^2,\qquad x\in L_p(\M).
\]
For $2<p<\infty$, the inequality is reversed.
\end{cor}

Another possible iteration is the following:

%co4 #&#
\begin{cor}
Let\vspace*{1pt} $(\M_n)_{1\le n\le N}$ be a family of von Neumann subalgebras of~$\M$. Assume that each $\M_n$ is $\s^\phi$-invariant and
$\phi\mid_{\M_n}$ is semifinite. Let $\E_n^+$ be the conditional
expectation with respect to $\M_n$ and $\E_n^{-}={\mathrm{Id}}-\E_n^+$.
Then for \mbox{$1<p\le2$},
\[
\llVert x\rrVert_p^2\ge\sum
_{(\varepsilon_i)\in\{+,-\}^N} (p-1)^{\llvert \{i\mid
\varepsilon_i=-\}\rrvert } \Biggl\llVert\Biggl(\prod
_{i=1}^N \E_i^{\varepsilon_i} \Biggr)
(x)\Biggr\rrVert_p^2,\qquad x\in L_p(\M).
\]
For $2<p<\infty$, the inequality is reversed.
\end{cor}

Applying it to the case where $\M=L_\infty(\{\pm1\}^N)$ and $\M_n$
is the subalgebra of functions
independent of the $n$th variable, we deduce the classical
hypercontractivity for the Walsh system (with operator valued
coefficients). Similarly, taking $\M$ to be the Clifford algebra with $N$
generators, we obtain the optimal hypercontractivity for Fermi
fields as pointed out in \cite{bcl,cl}.

We end the paper with some applications to hypercontractivity for group
von Neumann algebras.
In particular for the Poisson semigroup of a free group,
we obtain the optimal time for the hypercontractivity from $L_2$ to
$L_q$ for $q\ge4$.

%%%%%%%%%%%%%%%%%%%%%%%%%%%%%%%%%%%%%%%%%%%%%%%%%%%%%%
%%%%%%%%%%%%%%%%%%%%%%%%%%%%%%%%%%%%%%%%%%%%%%%%%%%%%%

%s2 #&#
\section{The proofs}

%%%%%%%%%%%%%%%%%%%%%%%%%%%%%%%%%%%%%%%%%%%%%%%%%%%%%%
%%%%%%%%%%%%%%%%%%%%%%%%%%%%%%%%%%%%%%%%%%%%%%%%%%%%%%

We will prove Theorems~\ref{conv-th} and \ref{BCL-th}. Using the
Haagerup reduction theorem as in \cite{hjx},
one can reduce both theorems to the finite case. Thus throughout this
section $\M$ will denote a von Neumann algebra
equipped with a faithful tracial normal state $\tau$. $L_p(\M)$ is
then constructed with respect to $\tau$.
We will first prove (\ref{BCL}), then deduce (\ref{conv}) from it.
$1<p<2$ will be fixed in the sequel.

As explained before, the proof of (\ref{BCL}) will be done by a
pseudo-differentiation argument.
Recall that for a continuous function
$f$ from an interval $I$ to $\real$ its pseudo-derivative of second
order at $t$ is
\[
D^2f(t)=\liminf_{h\to0^+} \frac{f(t+h)+f(t-h)-2f(t)}{h^2}.
\]
This pseudo-derivative shares many properties of the second derivative.
For instance,
if $D^2f$ is nonnegative on $I$, then
$f$ is convex. Indeed, by adding $\varepsilon t^2$ to $f$ (with
$\varepsilon>0$), we can
assume that $D^2f(t)$ is positive for all $t$. If $f$ was not convex,
there would exist $ t_0<t_1$ in $I$ such that the function $f-g$
takes a positive value at some point of $(t_0, t_1)$, where $g$ is
the straight line joining the two points $(t_0, f(t_0))$ and $(t_1,
f(t_1))$. So $f-g$ achieves a local maximum at a point $s\in(t_0,
t_1)$. Consequently, $D^2f(s)=D^2(f-g)(s)\le0$, which is a
contradiction.

Our pseudo-differentiation argument consists in proving the following
inequality for $x,y\in L_p(\M)$:
{\renewcommand{\theequation}{$D^2_{x, y}$}
\begin{equation}
D^2\llVert x+ty\rrVert^2_p(0)\geq2(p-1)
\llVert y\rrVert_p^2.
\end{equation}}\setcounter{equation}{3}%
Here the differentiation is, of course, taken with respect to the
variable $t$. The arguments from \cite{bcl} can be adapted to give:

%le5 #&#
\begin{lem}\label{step1}
Let $a, b\in\M$ be self-adjoint elements with $a$ invertible. Then
$(D^2_{a, b})$ holds.
\end{lem}

\begin{pf}
As $a$ is invertible in $\M$, $a+tb$ is also invertible for
small $t$. Introduce an
auxiliary function $\psi$ on $\real$,
\[
\psi(t)=\llVert a+tb\rrVert_p^p=\t\bigl(
\bigl(a^2+t(ab+ba)+t^2b^2\bigr)^{p/2}
\bigr).
\]
$\psi$ is differentiable in a neighborhood of the origin and
\[
\psi'(t)=\frac{p}2\t\bigl[\bigl(a^2+t(ab+ba)+t^2b^2
\bigr)^{p/2-1} \bigl((ab+ba)+2tb^2\bigr) \bigr].
\]
As in \cite{bcl} by functional calculus, the operator
$(a^2+t(ab+ba)+t^2b^2)^{p/2-1}$ admits the following integral
representation:
%
%e4 #&#
\begin{eqnarray}\label{int-rep}
&& \bigl(a^2+t(ab+ba)+t^2b^2
\bigr)^{p/2-1}
\nonumber\\[-8pt]\\[-8pt]\nonumber
&&\qquad =c_p\int_0^\infty
s^{p/2-1}\frac{1}{s+a^2+t(ab+ba)+t^2b^2} \,ds,
\end{eqnarray}
where
\[
c_p^{-1}=\int_0^\infty
s^{p/2-1}\frac{1}{s+1} \,ds.
\]
Thus $\psi$ is twice differentiable at $t=0$ and
%
%e5 #&#
\begin{eqnarray}\label{int1}
\psi''(0) &=& p \t\bigl(\llvert a\rrvert
^{p-2}b^2 \bigr)
\nonumber\\[-8pt]\\[-8pt]\nonumber
&&{}  -c_p\int_0^\infty
s^{p/2-1}\t\biggl[\frac{1}{s+a^2} (ab+ba) \frac{1}{s+a^2} (ab+ba)
\biggr]\,ds.
\end{eqnarray}
It then follows that $\f=\psi^{2/p}$ is also twice
differentiable at $t=0$ and
\[
\f''(0)=\frac{2}p \biggl(\frac{2}p-1
\biggr)\llVert a\rrVert_p^{2-2p}\psi'(0)^2+
\frac{2}p \llVert a\rrVert_p^{2-p}
\psi''(0)\ge\frac{2}p \llVert a\rrVert
_p^{2-p}\psi''(0).
\]
Hence $(D^2_{a, b})$ will be a consequence of
%
%e6 #&#
\begin{equation}
\label{df1} \frac{1}p \llVert a\rrVert_p^{2-p}
\psi''(0)\ge(p-1)\llVert b\rrVert_p^2.
\end{equation}
To prove the last inequality we claim that $\psi''(0)$ increases when
$a$ is replaced by~$\llvert a\rrvert $.
Indeed, the trace inside the integral in (\ref{int1}) is equal to
twice the following sum:
\[
\t\biggl[\frac{a}{s+ a^2} b \frac{a}{s+ a^2} b \biggr] +\t\biggl[
\frac{a^2}{s+ a^2} b \frac{1}{s+a^2} b \biggr].
\]
The second term above depends only on $\llvert a\rrvert $ (recalling
that $a$ is
self-adjoint). It remains to
show that the first one increases when $a$ is replaced by $\llvert
a\rrvert $. By
decomposing $a$ into
its positive and negative parts, we see that the first term is equal to
\[
\t\biggl[\frac{a_+}{s+ a_+^2} b \frac{a_+}{s+ a_+^2} b \biggr] + \t
\biggl[
\frac{a_-}{s+ a_-^2} b \frac{a_-}{s+ a_-^2} b \biggr] -2\t\biggl[\frac
{a_+}{s+ a_+^2} b
\frac{a_-}{s+ a_-^2} b \biggr].
\]
All above traces are nonnegative. Therefore, the above quantity
increases when the subtraction is replaced
by addition. Then tracing back the argument and noting that
$\llvert a\rrvert =a_++a_-$, we get the desired inequality
\[
\t\biggl[\frac{a}{s+ a^2} b \frac{a}{s+ a^2} b \biggr] \le\t\biggl[
\frac{\llvert a\rrvert }{s+ a^2} b \frac{\llvert a\rrvert }{s+ a^2}
b \biggr].
\]
Returning back to (\ref{int1}), we deduce the claim. Thus in the
following we will assume that $a$ is a positive invertible element of
$\M$.

The positivity of $a$ will facilitate the calculation of $\psi''(0)$
as explained in \cite{bcl}. Since $a+tb$ is positive for small $t$, we have
\[
\psi(t)=\t\bigl((a+tb)^p \bigr).
\]
Thus for $t$ close to 0,
\[
\psi'(t)=p \t\bigl((a+tb)^{p-1}b \bigr).
\]
To calculate the second derivative we use the following integral representation:
\[
(a+tb)^{p-1}=d_p\int_0^\infty
s^{p-1} \biggl[\frac{1}s - \frac{1}{s+a+tb} \biggr]\,ds.
\]
Consequently,
\[
\psi''(0) =p\, d_p\int
_0^\infty s^{p-1}\t\biggl[
\frac{1}{s+a} b \frac{1}{s+a} b \biggr]\,ds.
\]
As shown in \cite{bcl}, the function
\[
F\dvtx  z\mapsto\t\biggl[\frac{1}{s+z} b \frac{1}{s+z} b \biggr]
\]
is convex on the positive cone of $\M$.

Let $u$ be the unitary operator in the polar decomposition of $b$ (as
$\M$ is finite, the usual partial isometry in this decomposition can
be chosen to be a self-adjoint unitary). Then clearly
\[
F(z)= \frac{1}2 \bigl( F(z) + F(uzu) \bigr)\geq F \biggl(
\frac{z+uzu}2 \biggr).
\]
Now $z'=\frac{z+uzu} 2$ commutes with $u$, so
\[
F(z)\geq F\bigl(z'\bigr)=\t\biggl[\frac{1}{s+z'} \llvert b
\rrvert\frac{1}{s+z'} \llvert b\rrvert\biggr].
\]
Let $\B$ be the Abelian von Neumann subalgebra of $\M$ generated by $b$,
and let $\E_b$ be the associated trace preserving conditional
expectation. Then
\[
F\bigl(z'\bigr)=\t\biggl[\E_b \biggl(
\frac{1}{s+z'} \llvert b\rrvert\frac{1}{s+z'} \biggr) \llvert b\rrvert
\biggr].
\]
However, the Kadison--Schwarz inequality implies
\[
\E_b \biggl(\frac{1}{s+z'} \llvert b\rrvert\frac{1}{s+z'}
\biggr)\ge\E_b \biggl(\frac{1}{s+z'} \llvert b\rrvert
^{1/2} \biggr) \E_b \biggl(\llvert b\rrvert
^{1/2} \frac{1}{s+z'} \biggr) =\E_b \biggl(
\frac{1}{s+z'} \biggr)^2\llvert b\rrvert.
\]
Hence, by the positivity of the trace on products of positive elements,
we deduce
\[
F(z)\geq\t\biggl[\E_b \biggl(\frac{1}{s+z'}
\biggr)^2\llvert b\rrvert^2 \biggr].
\]
Then by the operator convexity of $\frac{1} t$, we have
\[
\E_b \biggl(\frac{1}{s+z'} \biggr)\geq\frac{1}{s+\E_b(z')}=
\frac{1}{s+\E_b(z)}.
\]
Letting $\widetilde a=\E_b(a)$, we have just shown
%
%e7 #&#
\begin{equation}
\label{df2} \psi''(0) \ge p\,d_p\int
_0^\infty s^{p-1}\t\biggl[
\frac{1}{s+\widetilde a} b \frac{1}{s+\widetilde a} b \biggr]\,ds =
\widetilde
\psi''(0),
\end{equation}
where
\[
\widetilde\psi(t)=\t\bigl((\widetilde a+tb)^p \bigr).
\]
Finally, (\ref{df1}) immediately follows from (\ref{df2}). Indeed, by
(\ref{df2}) and the H\"older inequality,
\[
\frac{1}p \llVert a\rrVert_p^{2-p}
\psi''(0)\ge\frac{1}p \llVert\widetilde a
\rrVert_p^{2-p}\widetilde\psi''(0)
=(p-1)\llVert\widetilde a\rrVert_p^{2-p}\t\bigl(\widetilde
a^{p-2}b^2\bigr)\ge(p-1)\llVert b\rrVert
_p^2.
\]
This finishes the proof of the lemma.
\end{pf}%\hfill$\Box$

For $a\in\M$ self-adjoint, we denote by $s(a)=\un_{(0,\infty
)}(\llvert a\rrvert )$. $s(a)$ is the support of~$a$, that is,
the least projection $e$ of $\M$ such that $ea=a$. We say that $a$ has
full support if $s(a)=1$.

%le6 #&#
\begin{lem}\label{step2}
Let $a, b\in\M$ be self-adjoint with $s(a)=1$. Then $(D^2_{a, b})$ holds.
\end{lem}

\begin{pf}
We will reduce this lemma to the previous one by cutting $a+tb$
with the spectral projections of $a$.
Let $e$ be a nonzero spectral projection of $a$, and put $a_e=eae$ and
$b_e=ebe$. Since $a$ is of full support, $a_e$ is invertible in the
reduced von Neumann algebra $\M_e=e\M e$. Thus Lemma \ref{step1}
can be applied to the couple $(a_e, b_e)$ in $\M_e$. Let
$\psi_e(t)=\llVert a_e+tb_e\rrVert _p^2$ as before. $\psi_e$ is twice
differentiable at $t=0$,
and (\ref{df1}) holds with $\psi_e$ in place of $\psi$.

Let $e^\perp=1-e$. Then for $t$ in a neighborhood of the origin, we
have [recalling that $\f(t)=\llVert a+tb\rrVert _p^2$]
\[
\f(t)\ge\bigl(\bigl\llVert e(a+tb)e\bigr\rrVert_p^p+
\bigl\llVert e^\perp(a+tb)e^\perp\bigr\rrVert
_p^p \bigr)^{2/p} {\mathop=^{\mathrm{def}}}
\bigl(\psi_e(t)+\g_e(t) \bigr)^{2/p}.
\]
However,
\[
\psi_e(t)=\llVert a_e\rrVert_p^p+t
\psi_e'(0)+\frac{t^2}2 \psi_e''(0)+{\mathrm{o}}\bigl(t^2\bigr)\qquad\mbox{as } t\to0.
\]
Let
\[
\alpha(t)=\llVert a_e\rrVert_p^p+
\g_e(t)=\bigl\llVert a_e+e^\perp(a+tb)e^\perp
\bigr\rrVert_p^p=\bigl\llVert a+te^\perp
be^\perp\bigr\rrVert_p^p.
\]
Then
\begin{eqnarray*}
\f(t) &\ge&\biggl(\alpha(t)+t\psi_e'(0)+
\frac{t^2}2 \psi_e''(0)+{\mathrm{o}}
\bigl(t^2\bigr) \biggr)^{2/p}
\\
&=&\alpha(t)^{2/p} \biggl(1+\frac{2t}p \frac{\psi_e'(0)}{\alpha
(t)}+\frac{t^2}{p} \frac{\psi_e''(0)}{\alpha(t)} +\frac{1}p \biggl(
\frac{2}p-1 \biggr)t^2\frac{\psi_e
'(0)^2}{\alpha(t)^2}+{\mathrm{o}}
\bigl(t^2\bigr) \biggr)
\\
&=&\alpha(t)^{2/p}+\frac{2t}p \psi_e'(0)
\alpha(t)^{2/p-1}+\frac{t^2}{p} \psi_e''(0)
\alpha(t)^{2/p-1}
\\
&&{} +\frac{1}p \biggl(\frac{2}p-1
\biggr)t^2\psi_e'(0)^2\alpha
(t)^{2/p-2}+{\mathrm{o}}\bigl(t^2\bigr)
\\
&\ge&\alpha(t)^{2/p}+\frac{2t}p \psi_e'(0)
\alpha(t)^{2/p-1}+\frac{t^2}{p} \psi_e''(0)
\alpha(t)^{2/p-1}+{\mathrm{o}}\bigl(t^2\bigr).
\end{eqnarray*}
By convexity of norms,
\[
\alpha(t)^{2/p}+\alpha(-t)^{2/p}\ge2\llVert a\rrVert
_p^2=2\f(0).
\]
We then deduce that
\begin{eqnarray*}
&& \frac{ \f(t)+\f(-t)-2\f(0)}{t^2}
\\
&&\qquad \ge\frac{2}p \psi_e'(0)
\frac{\alpha(t)^{2/p-1}-\alpha
(-t)^{2/p-1}}{t}
\\
&&\qquad\quad{}  +\frac{1}{p} \psi_e''(0)
\bigl[\alpha(t)^{2/p-1}+\alpha(-t)^{2/p-1} \bigr] +{\mathrm{o}}(1).
\end{eqnarray*}
The uniform smoothness of the norm $\llVert \rrVert _p$ implies that
the function
$\alpha^{2/p-1}$ is differentiable
at $t=0$, and its derivative is equal to
\[
(2-p)\llVert a\rrVert_p^{1-p}\t\bigl(v\llvert a\rrvert
^{p-1}e^\perp be^\perp\bigr) {\mathop=^\mathrm{def}}
\d_e,
\]
where $v$ is the unitary in the polar decomposition of $a$. It then
follows that
\[
D^2\f(0)\ge\frac{4}p \psi_e'(0)
\d_e+\frac{2}{p} \psi_e''(0)
\llVert a\rrVert_p^{2-p}.
\]
Hence by (\ref{df1}),
\[
D^2\f(0)\ge\frac{4}p \psi_e'(0)
\d_e+2(p-1)\llVert b_e\rrVert_p^2.
\]
Thanks to the full support assumption of $a$, we can let $e\to1$ in
the above inequality. This limit procedure
removes the first extra term, so we finally get
\[
D^2\f(0)\ge2(p-1)\llVert b\rrVert_p^2.
\]\upqed
\end{pf}
%
%\hfill$\Box$

Now we are ready to show (\ref{BCL}).

\begin{pf*}{Proof of Theorem~\ref{BCL-th}}
First by
density, we need only to show (\ref{BCL})
for \mbox{$x, y\in\M$}. Then notice that it suffices to do it for
self-adjoint elements using a classical $2\times2$-matrix trick.
Indeed, let $\widetilde\M=\ma_2\otimes\M$ equipped with the tensor
trace. Given $x, y\in\M$ let
\[
a= %
\pmatrix{ 0&x
\cr
x^*&0} \quad\mbox{and}\quad b= %
\pmatrix{ 0&y
\cr
y^*&0}.
\]
Then $a$ and $b$ are self-adjoint. Moreover, by easy computations, (\ref
{BCL}) for $x$ and $y$
is equivalent to the same inequality for $a$ and $b$.

To use Lemma \ref{step2}, we require that $a+tb$ have full
support for any $t\in\real$. This is achieved by a tensor
product argument. Choose a positive element $c\in L_\infty([0, 1])$ whose
spectral measure with respect to Lebesgue measure is diffuse
(atomless), say $c(t)=t$ for $t\in[0, 1]$. In other words, considered
as a random variable in the probability space $[0, 1]$, the law of
$c$ is diffuse. On the other hand, for $t\in\real$, composing the spectral
resolution of $a+tb$ with the trace $\tau$, we can view $a+tb$ as a
random variable in another probability space $(\O, P)$. Now, consider
the tensor von Neumann algebra $L_\infty([0, 1])\,\overline\otimes\,\M
$; it is finite.
$\M$ and
$L_\infty([0, 1])$ are identified as subalgebras of
$L_\infty([0, 1])\,\overline\otimes\,\M$ in the usual way. Then for
any $\varepsilon>0$,
the law of $a+tb+\varepsilon c$ is the convolution of the laws of $a+tb$
and $\varepsilon c$. It is atomless since the law of $c$
is atomless. Consequently, the support of $a+tb+\varepsilon c$ is full in
$L_\infty([0, 1])\,\overline\otimes\,\M$.\vadjust{\goodbreak}

Thus by Lemma~\ref{step2} applied to the pair $(a+\varepsilon c, b)$,
the function $f(t)= \llVert a+tb+\varepsilon
c\rrVert _p^2-(p-1)t^2\llVert b\rrVert _p^2$ satisfies $D^2(f)(t)\geq
0$, so it is
convex.\vspace*{1pt} Hence
\mbox{$f(1)+f(-1)\geq2 f(0)$;} this is (\ref{BCL}) for $a+\varepsilon c$ and
$b$. Letting
$\varepsilon\to0$ gives the desired result.
\end{pf*} %\hfill$\Box$

Finally, we deduce (\ref{conv}) from (\ref{BCL}).

\begin{pf*}{Proof of Theorem~\ref{conv-th}}
Given
$x\in L_p(\M)$ let $a=\E(x)$ and $b=x-\E(x)$.
Consider again the function $f$ defined by
\[
f(t)=\llVert a+tb\rrVert_p^2-(p-1)t^2
\llVert b\rrVert_p^2.
\]
Then (\ref{BCL}) implies $D^2f\ge0$, so $f$ is convex. On the other hand,
the function $g(t)=\llVert a+tb\rrVert _p^2$ is also convex and by the
contractivity of $\E$ on $L_p(\M)$,
\[
\llVert a+tb\rrVert_p\ge\bigl\llVert\E(a+tb)\bigr\rrVert
_p=\llVert a\rrVert_p.
\]
Hence we conclude that the right derivative $g_r'(0)\ge0$, so that
$f_r'(0)\ge0$, too. Consequently, $f$ is increasing on $\real^+$. In
particular, $f(1)\geq f(0)$, which is nothing but (\ref{conv}).
\end{pf*}
%\hfill$\Box$

%%%%%%%%%%%%%%%%%%%%%%%%%%%%%%%%%%%%%%%%%%%%%%%%%%%%%%
%%%%%%%%%%%%%%%%%%%%%%%%%%%%%%%%%%%%%%%%%%%%%%%%%%%%%%

%s3 #&#
\section{Applications to hypercontractivity}

%%%%%%%%%%%%%%%%%%%%%%%%%%%%%%%%%%%%%%%%%%%%%%%%%%%%%%
%%%%%%%%%%%%%%%%%%%%%%%%%%%%%%%%%%%%%%%%%%%%%%%%%%%%%%

We give in this section some applications to hypercontractivity
inequalities on group von Neumann algebras.
Let $G$ be a discrete group and $vN(G)$ the associated group von
Neumann algebra. Recall that $vN(G)$ is the von Neumann algebra
generated by the left regular representation $\l$:
$vN(G)=\l(G)''\subset B(\el_2(G))$. It is equipped with a canonical
trace $\t$, that is, $\t(x)=\langle xe, e\rangle$, where $e$ is the
identity of
$G$. Given a function $\psi\dvtx  G\to\real_+$ with $\psi(e)=0$, we consider
the associated Fourier--Schur multiplier initially defined on the
family $\com[G]$ of polynomials on $G$:
\[
P_t\dvtx  \sum_{g\in G}x(g)\l(g)\mapsto\sum
_{g\in G}e^{-t\psi(g)}x(g)\l(g),\qquad t>0.
\]
We will assume that $P_t$ extends to a contraction on $L_p(vN(G))$ for
every $1\le p\le\infty$.
Schoenberg's classical theorem asserts that if $\psi$ is symmetric and
conditionally negative,
$P_t$ is a completely positive map on $vN(G)$. Since it is trace
preserving, $P_t$ defines a contraction
on $L_p(vN(G))$ for every $1\le p\le\infty$. Thus in this case our
assumption is satisfied.

The
hypercontractivity problem for the semigroup $(P_t)_{t>0}$ and for $1 <
p <q < \infty$, consists in determining the optimal
time $t_{p,q} > 0$ such that
\[
\llVert P_t\rrVert_{p\to q}\le1\qquad\forall t\ge
t_{p,q}.
\]
We refer to
\cite{jppp,jpppr} for more information and historical references. It
is easy to check that if such a time $t_{p,q}$ exists, then $\psi$ has a
spectral gap, namely \mbox{$\inf_{g\in G\setminus\{e\}}\psi(g)>0$}. After
rescaling, we will assume that $\inf_{g\in G\setminus\{e\}}\psi(g)=1$.

In most-known cases the expected optimal time $t_{p,q}$ is attained, namely,
\[
t_{p,q}=\log\sqrt{\frac{q-1}{p-1}}.
\]
It is a particularly interesting problem of determining the optimal
time $t_{p, q}$ when $G=\FF_n$ is the free group on $n$ generators
with $n\in\nat\cup\{\infty\}$, and $\psi$ is its natural length function.
Some partial results are obtained in \cite{jppp,jpppr}. For instance,
by embedding $vN(\FF_n)$ into
a free product of Clifford algebras, it is proved in \cite{jpppr} that
for any $q>2$,\vspace*{-2pt}
\[
t_{2, q}\le\log\sqrt{q-1}+ \biggl(\frac{1}2-\frac{1}q
\biggr)\log\sqrt2.
\]
On the other hand, Junge et al. \cite{jppp} show that for any finite
$n$ there exists $q(n)$ such that if $q\ge q(n)$ is an even integer, then
\[
t_{2, q}= \log\sqrt{q-1}.
\]
The proof is combinatoric and based on lengthy calculations.

Here we provide an improvement. We will use Haagerup-type inequalities
\cite{haag}. Letting $S_k$ be the set of words of length $k$ in $\FF
_n$ and for any $x\in vN(\FF_n)$ supported
on $S_k$, the original Haagerup inequality is
%
%e8 #&#
\begin{equation}
\label{K8} \llVert x\rrVert_\infty\le(k+1)\llVert x\rrVert
_2.
\end{equation}
For $q>2$ and $k\in\nat$, let $K_{k,
q}$ be the best constant in the following Khintchine inequality for
homogeneous polynomials $x$ of degree $k$:
\[
\biggl\llVert\sum_{g\in S_k}x(g)\l(g)\biggr\rrVert
_q\le K_{k, q}\biggl\llVert\sum
_{g\in S_k}x(g)\l(g)\biggr\rrVert_2.
\]
We will need the following:

%le7 #&#
\begin{lem}\label{kin4}
We have $K_{k,4}\leq(k+1)^{1/4}$.
\end{lem}

\begin{pf}
Denote by $g_i$ the generators of $\FF_n$ with the convention that
$g_{-i}=g_i^{-1}$. For a multi-index
$\underline i =(i_1,\ldots,i_d)$ with $i_j+i_{j+1}\neq0$, we let
$g_{\underline i}= g_{i_1}\cdots g_{i_d}$ and $\llvert \underline
i\rrvert =d$. So we may write\vspace*{2pt}
\[
x=\sum_{\llvert \underline i\rrvert =k} \alpha_{\underline i}
\l(g_{\underline i}).
\]
We compute $x^*x$ according to simplifications that may occur
\[
x^*x= \sum_{0\leq d\leq k} \mathop{\mathop{\sum
_{\llvert \underline\beta\rrvert = d}}_{\llvert \underline i\rrvert
=\llvert \underline
j\rrvert =k-d}}_{i_{k-d}\neq j_{k-d}\neq-\beta_1} \overline{\alpha_{\underline
j,\underline\beta}}
\alpha_{\underline i,\underline\beta} \l\bigl(g_{\underline
j}^{-1}g_{\underline i}
\bigr),
\]
where $\underline i,\underline\beta$ denotes the multi-index obtained
by superposing the two
multi-indices $\underline i$ and $\underline\beta$.
We then deduce that
\begin{eqnarray*}
\llVert x\rrVert_4^4=\bigl\llVert x^*x\bigr\rrVert
_2^2 & = &\sum_{0\leq d\leq k}
\mathop{\sum_{\llvert \underline i\rrvert =\llvert \underline
j\rrvert =k-d }}_{i_{k-d}\neq j_{k-d}} \biggl(\mathop{
\sum_{\llvert \underline\beta\rrvert = d}}_{i_{k-d}\neq
j_{k-d}\neq-\beta_1} \overline{
\alpha_{\underline j,\underline
\beta}} \alpha_{\underline i,\underline\beta} \biggr)^2
\\
& \leq&\sum_{0\leq d\leq k} \mathop{\sum
_{\llvert \underline
i\rrvert =\llvert \underline j\rrvert =k-d}}_{i_{k-d}\neq j_{k-d}}
\biggl(\mathop{\sum
_{\llvert \underline\beta\rrvert = d}}_{
i_{k-d} \neq-\beta_1} \llvert\alpha_{\underline i,\underline\beta
}\rrvert
^2 \biggr)\cdot \biggl(\mathop{\sum_{\llvert \underline\beta\rrvert =
d}}_{j_{k-d} \neq-\beta_1}
\llvert\alpha_{\underline j,\underline\beta}\rrvert^2 \biggr)
\\
&\leq&\sum_{0\leq d\leq k} \biggl(\mathop{\sum
_{\llvert \underline
i\rrvert =k-d, \llvert \underline\beta\rrvert = d}}_{i_{k-d} \neq
-\beta_1} \llvert\alpha_{\underline i,\underline\beta}\rrvert
^2 \biggr)\cdot \biggl(\mathop{\sum_{\llvert \underline j\rrvert =k-d,
\llvert \underline\beta\rrvert = d}}_{j_{k-d}
\neq-\beta_1}
\llvert\alpha_{\underline j,\underline\beta}\rrvert^2 \biggr)
\\
&=& (k+1) \llVert x\rrVert_2^4.
\end{eqnarray*}\upqed
\end{pf}%\qed

%re8 #&#
\begin{remark}
Taking $\alpha_{\underline i}=1$ and by the free central limit theorem
as $n\to\infty$, one can see that the previous inequality
is sharp. Thus $K_{k,4}=(k+1)^{1/4}$. This constant is the
$L_4$-norm of the $k$th
Chebyshev polynomial for the semi-circle law.
\end{remark}

Using the H\"older inequality we deduce from (\ref{K8}) and the
previous lemma that for any $k\ge1$,
%
%e9 #&#
%e10 #&#
\begin{eqnarray}
K_{k, q}&\le& (k+1)^{1-3/q},\qquad q\ge4, \label{qsup4}
\\
K_{k, q}&\le& (k+1)^{1/2 -1/q},\qquad2\leq q\leq4.\label{qinf4}
\end{eqnarray}

We will also use the following elementary folklore:

%re9 #&#
\begin{remark}\label{Lp-Lq}
Let $T\dvtx  L_p(\M)\to L_q(\M)$ be a bounded linear map. Assume that $T$
is $2$-positive in the sense that ${\mathrm{Id}}_{\ma_2}\otimes T$
maps the positive cone of $L_p(\ma_2\otimes\M)$ to that of $L_q(\ma
_2\otimes\M)$. Then
\[
\bigl\llVert T(x)\bigr\rrVert_q\le\bigl\llVert T\bigl(\llvert x
\rrvert\bigr)\bigr\rrVert_q^{1/2} \bigl\llVert T\bigl(\bigl
\llvert x^*\bigr\rrvert\bigr)\bigr\rrVert_q^{1/2},\qquad x
\in L_p(\M).
\]
Consequently,
\[
\llVert T\rrVert=\sup\bigl\{\bigl\llVert T(x)\bigr\rrVert_q\dvtx  x\in
L_p(\M)^+, \llVert x\rrVert_p\le1\bigr\}.
\]
\end{remark}

Indeed, for any $x\in L_p(\M)$,
\[
\pmatrix{ \llvert x\rrvert&x
\cr
x^*&\bigl\llvert x^*\bigr\rrvert}
\geq0.
\]
So the $2$-positivity of $T$ implies
\[
\pmatrix{ T\bigl(\llvert x\rrvert\bigr)& T(x)
\vspace*{3pt}\cr
T\bigl( x^*\bigr)&
T\bigl(\bigl\llvert x^*\bigr\rrvert\bigr)}%
\geq0.
\]
This yields a contraction $c\in\M$ such that $T(x)=T(\llvert x\rrvert
)^{1/2}c
T(\llvert x^*\rrvert )^{1/2}$. Then the H\"older inequality
gives the assertion.

%th10 #&#
\begin{teo}\label{freehyper}
There exists $\varepsilon_0>0$ such that for any free $\FF_n$ and any
\mbox{$q\ge4-\varepsilon_0$},
\[
\llVert P_t\rrVert_{2\to q}\le1\quad\Leftrightarrow\quad t\ge\log{
\sqrt{q-1}}.
\]
\end{teo}

\begin{pf}
The necessity is clear. The proof of the sufficiency will rely on
Remark 3.7 of \cite{jpppr}. Let $\sigma$ be the
automorphism of $vN(\FF_n)$ given by $\sigma(\l(g_i))=\l(g_{i}^{-1})$.
Then $P_t$ is hypercontractive from $L_2$ to $L_q$ with optimal
time on $vN(\FF_n)^{\sigma}$, the fixed point algebra of $\sigma$.
Let $\E$ be the conditional expectation onto $vN(\FF_n)^{\sigma}$. Note
that $\E= \frac{{\mathrm{Id}} +\sigma}2$ and it commutes with $P_t$.

Fix $q>2$. To prove $\llVert P_t\rrVert _{2\to q}\le1$ for
$t\ge\log{\sqrt{q-1}}$, it suffices to show $\llVert P_t(x)\rrVert
_{q}\leq\llVert x\rrVert
_2$ for any positive $x\in\com[\FF_n]$
by virtue of Remark~\ref{Lp-Lq}.
We need one more reduction. Given complex numbers $\zeta_i$ of modulus
1, there exists an
automorphism $\pi_\zeta$ of $vN(\FF_n)$ given by $\pi(g_i)=\zeta_i g_i$.
It is an isometry on all $L_p$'s.
Note that $\pi_\zeta$ and $P_t$ commute. Thus to prove
$\llVert P_t(x)\rrVert _{q}\leq\llVert x\rrVert _2$, we may assume
that $x(g_i)$ is real for every generator $g_i$. We will fix a positive
$x\in\com[\FF_n]$ with the last property.

Then write $x=y+z$ where $y=\E(x)$.
Since $x(g_i)\in\real$, we have that $z$ does not have constant terms nor
of degree $1$. By Theorem \ref{conv-th} (or Theorem \ref{BCL}) and
Remark~3.7 of \cite{jpppr},
\[
\bigl\llVert P_t(x)\bigr\rrVert_q^2 \leq
\bigl\llVert P_t(y)\bigr\rrVert_q^2 + (q-1)
\bigl\llVert P_t(z)\bigr\rrVert_q^2 \leq
\llVert y\rrVert_2^2 + (q-1)\bigl\llVert
P_t(z)\bigr\rrVert_q^2.
\]
Then for $t=\log\sqrt{q-1}$, decomposing $z$ according to its
homogeneous components $(z_k)$ and using the Khintchine and the
Cauchy--Schwarz inequalities, we get
\[
\bigl\llVert P_t(z)\bigr\rrVert_q^2 \le
\biggl(\sum_{k\ge2}e^{-tk} \llVert
z_k\rrVert_q \biggr)^2 \le\sum
_{k\ge2}K_{k, q}^2 \frac{1}{(q-1)^k}
\llVert z\rrVert_2^2.
\]
We aim to find those $q>2$ for which
\[
R_q=(q-1)\sum_{k\ge2}K_{k, q}^2
\frac{1}{(q-1)^k}\le1.
\]
For $q\ge4$, by (\ref{qsup4}) we have
\[
R_q \le\sum_{k\ge2}(k+1)^{2(1-3/q)}
\frac{1}{(q-1)^{k-1}}.
\]
The terms of the sum on the right-hand side are decreasing functions of
$q$ if their derivatives are negative, that is, if
\[
\frac{6(q-1)}{q^2}\le\frac{k-1}{\log{(k+1)}}.
\]
Noting that the
left-hand side of the above inequality is decreasing on $q$, one
easily checks that this inequality is true for $q\ge4$ and
$k\ge3$. However, it is true for $k=2$ if and only if $q\ge q_0$, where
\[
q_0=\sqrt{3\log3} (\sqrt{3\log3}+\sqrt{3\log3 -2})\approx5.36244.
\]
We have the following numerical estimates:
\[
R_4\leq0.92952\quad\mbox{and}\quad\frac{3^{2(1-3/{q_0})}}{q_0-1}-
\frac{3^{1/2}}3\leq0.02613.
\]
Hence if $q\in[4,q_0]$,
\[
R_{q}\leq R_4 +\frac{3^{2(1-3/{q_0})}}{q_0-1}- \frac{3^{1/2}}3<1.
\]
We thus conclude that $R_q<1$ for all $q\geq4$.

Since $R_q$ is dominated by a continuous function of $q$, using (\ref
{qinf4}) we get a similar estimate for
$q\geq4-\varepsilon_0$ for some $\varepsilon_0$.
A numerical estimate gives $\varepsilon_0\approx0.18$.
\end{pf}
%\qed

%re11 #&#
\begin{remark}
Instead of Remark 3.7 of \cite{jpppr}, we can equally use Theorem~A(iii) of \cite{jpppr} in the preceding proof. But the commutation
of $P_t$ and the conditional expectation onto the symmetric subalgebra
$\mathcal A_{\mathrm{sym}}^n$ is less obvious.
\end{remark}

It is likely that $\varepsilon_0=2$, but other methods would have to be
developed.

Gross's pioneering work \cite{g} shows that hypercontractivity is
equivalent to the validity of log-Sobolev inequalities.
In the present situation of free groups, the validity of the
hypercontractivity with optimal time in full generality
(or equivalently, $\varepsilon_0=2$) is equivalent to the following
log-Sobolev inequality in $L_q$ for any $q\ge2$:
{\renewcommand{\theequation}{$\mathrm{SL}_q$}
\begin{equation}\label{eqSLq}
\t\bigl(x^q\log x\bigr)\le\frac{q}{2(q-1)} \t\bigl(x^{q-1}L(x)
\bigr)+\llVert x\rrVert_q^q\log\llVert x\rrVert
_q, \qquad x\in\mathcal D^+.
\end{equation}}\setcounter{equation}{3}%
Here $L$ denotes the negative generator of $(P_t)$, and $\mathcal D$ is
a core for $L$ where the inequality makes sense. It is known
that $(\mathrm{SL}_2)$ implies $(\mathrm{SL}_q)$ for all $q$; see
\cite{oz}. In the same spirit we can show that $(\mathrm{SL}_p)$
implies $(\mathrm{SL}_q)$ if $q>p\ge2$. Let
us record this explicitly here since it might be of
interest. The semigroup $(P_t)$ can be any completely positive
symmetric Markovian semigroup such that $\mathcal D$ is rich enough.

%re12 #&#
\begin{remark}
Let $q>p\ge2$. Then $(\mathrm{SL}_p)$ implies $(\mathrm{SL}_q)$.
\end{remark}

To check the remark we rewrite $(\mathrm{SL}_q)$ in a symmetric form with
respect to $q$ and its conjugate index $q'$ (provided that $\mathcal
D$ is big enough):
{\renewcommand{\theequation}{$\mathrm{SL}^s_q$}
\begin{equation}
\t(x\log x)\le\tfrac{1}{2} q'q \t\bigl(x^{1/{q'}}L
\bigl(x^{1/q}\bigr) \bigr)+\tau(x)\log\tau(x),\qquad x\in\mathcal D^+.
\end{equation}}\setcounter{equation}{10}%
Recall that for $y\in\mathrm{Dom}(L)$, $\tau(zL(y) )=\lim_{r\to
0} \frac{1} r \t(z(1-P_r)(y) )$. Let $r>0$ and $x\in\M^+$,
and we will check that the function
$q\mapsto q'q \t(x^{1/{q'}}(1-P_r) (x^{1/q} )
)$ is
increasing for $q\ge2$; we put $\theta=\frac{1}q$. It is known from
\cite{ahk} that there exists a positive symmetric
Borel measure $\mu_r$ on $\s(x)\times\s(x)$ such that
\[
\t\bigl(x^{1-\theta}P_r\bigl(x^\theta\bigr) \bigr) =
\int_{\s(x)\times\s(x)}s^{1-\theta}t^\theta d
\mu_r(s, t).
\]
Hence, by symmetry, it suffices to show that
\[
f\dvtx \theta\mapsto\frac{1+ u-u^{\theta}-u^{1-\theta}}{\theta
(1-\theta)}
\]
is convex on $[0, 1]$ for $u>0$ as $f(\theta)=f(1-\theta)$. One
easily checks that
\begin{eqnarray*}
f(\theta) &=& \int_0^1 \log(u)
\bigl(u^{\theta+(1-\theta)(1-t)}- u^{\theta t}+u^{1-\theta+\theta
t}-u^{(1-\theta)(1-t)} \bigr)
\,dt,
\\
f''(\theta) &=& \int_0^1
\log(u)^3 \bigl(t^2 \bigl(u^{\theta+(1-\theta)(1-t)}-
u^{\theta
t} \bigr)
\\
&&\hspace*{51pt}{} +(1-t)^2 \bigl(u^{1-\theta+\theta
t}-u^{(1-\theta)(1-t)}
\bigr) \bigr) \,dt
\geq 0.
\end{eqnarray*}
Passing to the limit
in $r$ gives the result if $\mathcal D$ is big enough.

We end this section with application to more general groups $(G,\psi)$.
If $\psi$ is symmetric and satisfies the exponential order growth
%
%e11 #&#
\begin{equation}
\label{exp} \bigl\llvert\bigl\{g\in G\dvtx  \psi(g)\le R\bigr\}\bigr\rrvert
\le C
\rho^R\qquad\forall R>0
\end{equation}
for some $C>0$ and $\rho>1$, then one of the main results of \cite
{jppp} shows that for $2<q<\infty$,
\[
t_{2,q}\le\eta\log\sqrt{q-1}
\]
for any $\eta>2$ when $\rho$ is large compared to $C$. Their argument
consists in first considering the case $q=4$
by combinatoric methods and then using Gross's extrapolation. We will
show that the martingale inequality
in Theorem~\ref{conv-th} easily implies a slight improvement. Note
that our estimate on $t_{2,q}$ is
as close as to the expected optimal time as when $q$ is sufficiently
large, compared to $\rho$ and $C$.

%pr13 #&#
\begin{prop}
Assume (\ref{exp}) and $2<q<\infty$.
Then
\[
t_{2,q}\leq\biggl(\frac{q-2}q \log\sqrt{2C\rho}+\log\sqrt{q-1}
\biggr) \vee\log\rho.
\]
\end{prop}

\begin{pf}
By (\ref{exp}), the range of $\psi$ is countable. Let $\psi(G)=\{
n_0, n_1, n_2,\ldots\}$ with $n_0<n_1< n_2<\cdots.$
Then $n_0=0$ and $n_1=1$. Let $x\in vN(G)$ be a polynomial,
$x=\sum x(g)\l(g)$,
and let $y=x-x(e)$. By Theorem~\ref{conv-th}
\[
\bigl\llVert P_t(x)\bigr\rrVert_q^2\le
\bigl\llvert x(e)\bigr\rrvert^2+(q-1)\bigl\llVert P_t(y)
\bigr\rrVert_q^2.
\]
Let $B_k=\{g\in G\dvtx  \psi(g)\le n_k\}$, $S_k=B_k\setminus B_{k-1}$
and $y_k=\sum_{g\in S_k} x(g) \l(g)$.
Then
\[
\bigl\llVert P_t(y)\bigr\rrVert_\infty^2 \leq
\biggl(\sum_{k\ge1}e^{-tn_k} \llVert
y_k\rrVert_\infty\biggr)^2 \leq\biggl(\sum
_{k\ge1}e^{-2tn_k} \llvert S_k\rrvert
\biggr)\cdot \biggl(\sum_{k\ge
1}\frac{\llVert y_k\rrVert _\infty^2}{\llvert S_k\rrvert }
\biggr).
\]
Obviously,
\[
\llVert y_k\rrVert_\infty^2\leq\biggl(\sum
_{g\in S_k} \bigl\llvert x(g)\bigr\rrvert
\biggr)^2 \leq\llvert S_k\rrvert\sum
_{g\in S_k} \bigl\llvert x(g)\bigr\rrvert^2.
\]
We get, using the H\"older inequality,
\[
\bigl\llVert P_t(y)\bigr\rrVert_q^2 \leq
e^{-{4t}/q} \biggl(\sum_{k\ge1}e^{-2tn_k}
\llvert S_k\rrvert\biggr)^{(q-2)/q} \llVert y\rrVert
_2^2.
\]
Actually exchanging the arguments, one has the following, slightly
better estimate that we will not use:
\[
\bigl\llVert P_t(y)\bigr\rrVert_q^2 \leq
\sum_{k\ge1}e^{-2tn_k} \llvert S_k
\rrvert^{{2(q-2)}/q} \llVert y\rrVert_2^2.
\]
By (\ref{exp}), for $t>\log\rho$,
\begin{eqnarray*}
\sum_{k\ge1}e^{-2tn_k}\llvert S_k
\rrvert&\leq&\sum_{k\ge1} \bigl(e^{-2tn_k}-e^{-2tn_{k+1}}
\bigr)\llvert B_k\rrvert\le2C t\int_1^\infty
e^{-(2t-\log\rho)s} \,ds
\\
&=&2C \frac{t}{2t-\log\rho} e^{-(2t-\log\rho)} \le2C e^{-(2t-\log\rho)}.
\end{eqnarray*}
Hence,\vspace*{1pt} if $2t\ge\frac{q-2}{q}\log(2C\rho)+\log(q-1)$, we deduce
$\llVert P_t(x)\rrVert _q\le\llVert x\rrVert _2$,
whence the assertion.
\end{pf}
%\hfill$\Box$

%%%%%%%%%%%%%%%%%%%%%%%%%%%%%%%%%%%%%%%%%%%%%%%%%%%%%%
%%%%%%%%%%%%%%%%%%%%%%%%%%%%%%%%%%%%%%%%%%%%%%%%%%%%%%

%\begin{appendix}
%\section{}
%\end{appendix}

%\section*{Acknowledgments}

%\begin{supplement}[id=suppA]
%\slink[doi]{10.1214/00-AOPXXXXSUPP} %[doi,text={...}] - jei reikia
%suskaldyti doi
%\sfilename{aopXXXX\_supp.pdf}
%\end{supplement}

% imsref loaded by linak, 2014-12-23 11:58:05

\printaddresses
\end{document}